\documentclass[11pt]{amsart}
\usepackage{amsmath,amssymb,latexsym,dsfont}
\usepackage[left=2.6cm,right=2.6cm,top=2.7cm,bottom=2.7cm]{geometry}

    \usepackage{amsmath,amsfonts,amssymb,epsfig}
     \usepackage{graphicx,dsfont}
     \usepackage{hyperref}
     \usepackage{cases}
     \usepackage{color}

    \newcommand{\EE}{\mathcal{E}}  
          \newcommand{\HH}{\mathcal{H}}

      \newcommand{\be}{{\bf  e}}

\newtheorem{theorem}{Theorem}[section]
\newtheorem{lemma}{Lemma}[section]

\newtheorem{proposition}{Proposition}[section]

\newtheorem{remark}{Remark}[section]

\newtheorem{definition}{Definition}[section]
\numberwithin{equation}{section}

      \newcommand{\R}{{\mathbb{R}}}

      \newcommand{\curl}{\operatorname{curl}}
      \newcommand{\dive}{\operatorname{div}}
      
      \newcommand{\ext}{\operatorname{ext}}
      \newcommand{\inte}{\operatorname{int}}
      \newcommand{\loc}{\operatorname{loc}}
      
      \newcommand{\eps}{\varepsilon}
      \newcommand{\mR}{\mathbb{R}}
      
      \newcommand{\pbE}{\operatorname{\mathcal E}}
      \newcommand{\pbH}{\operatorname{\mathcal H}}

      \newcommand{\bE}{\operatorname{\textbf{E}}}
      \newcommand{\bH}{\operatorname{\textbf{H}}}

     \newcommand{\supp}{\mbox{supp }}
           \newcommand{\dsp}{\displaystyle}

      \makeatletter
      \def\@setcopyright{}
      \def\serieslogo@{}
      \makeatother

\begin{document}

    \author[H.-M. Nguyen]{Hoai-Minh Nguyen}
\author[M.S. Vogelius]{Michael S. Vogelius}

\address[H.-M. Nguyen]{Department of Mathematics, EPFL SB CAMA, Station 8,  \newline\indent
	 CH-1015 Lausanne, Switzerland.}
\email{hoai-minh.nguyen@epfl.ch}

\address[M.S. Vogelius]{Department of Mathematics, Rutgers University, 
	\newline\indent New Brunswick, NJ 08903, USA.}
\email{vogelius@math.rutgers.edu}
 \thanks{Michael S. Vogelius is partially supported by NSF grant DMS-12-11330.}
 \title[Approximate cloaking using transformation optics]{Approximate cloaking using transformation optics for acoustic and electromagnetic waves \\
}
%
%
%
%
%
%
%
%

  \maketitle
%


\tableofcontents

\section{Introduction}

 
Cloaking, using transformation optics, was introduced by Pendry, Schurig, and Smith \cite{Pen} for the Helmholtz equation, and by Leonhardt \cite{Leo} in the geometric optics setting. They employed a singular change of variables which blows up a point to form the
cloaked region, and correspondingly transforms the surrounding medium. The same transformation had been used  by Greenleaf, Lassas,
and Uhlmann  to establish (singular) non-uniqueness in Calderon's problem \cite{Green}.  Transformation cloaking essentially relies on the invariance of the Maxwell equations or the wave equation in the time harmonic as well as in the real time regime.  The singular, and very anisotropic, nature of the cloaks presents various difficulties in practice as well as in theory: (1) they are hard to fabricate, and (2) in certain cases the correct definition (and therefore the properties) of the corresponding acoustic/electromagnetic fields is an issue. To avoid using
the singular structure, various regularization schemes have been proposed. One of these  was suggested by Kohn, Shen, Vogelius, and Weinstein  in \cite{Kohn}, where they  used a  transformation which maps a small ball, instead of a point, to the cloaked region. For other, related regularization schemes see \cite{RYNQ,GKLU07}. 

In this paper we discuss various results on approximate cloaking for acoustic and electromagnetic waves, using schemes in the spirit of the one in \cite{Kohn}. We consider the time harmonic regime as well as the real time regime.   The discussion reflects our taste, our understanding,  and our works in this direction. We have not made any attempt to perform an extensive review of work in this active area.

Let us briefly describe  the structure and the contents of this paper. The first section after this introduction is devoted to approximate cloaking for acoustic waves in the time harmonic regime. In this section, we initially discuss approximate cloaking when employing a (damping) lossy layer, and  emphasize the role of frequency in the estimates for the degree of (in)visibility. We then discuss several results relating to the case where no lossy layer is present. The notion of resonance appears naturally here, and with it the possibility of a situation in which the energy inside the cloaked region goes to infinity and cloaking might not be achieved.  
Finally, we discuss various ways to enhance the cloaking effect. The results in this section are based on  \cite{Ng-1, Ng-Vogelius-2, Ng-2, GV, HV}. 
The following section concerns approximate cloaking for acoustic waves in the real time regime.  Here  we present a natural cloaking scheme   and discuss the general approach we have used  to estimate the degree of (in)visibility in the time regime. The central idea is to connect the problem in the time regime with its corresponding problem in the time harmonic regime, via the Fourier transform in time, and then use the analysis discussed in the previous section. A technical point of this approach is to establish the outgoing radiation condition for the Fourier transform (with respect to time) of solutions to the wave equation. The resolution of this technical point is interesting in itself and has been applied in different contexts, see, e.g., \cite{MinhLinh}.  The materials in this section are based on \cite{Ng-Vogelius-3, Ng-Vogelius-4}. 
The next section is on approximate cloaking for electromagnetic waves in the time harmonic regime. We immediately discuss  the situation in which no lossy layer is present.  The notion of resonance also naturally occurs 
in this case. However, in  contrast to the acoustic case,  cloaking is always achieved even in the presence of resonance. That being said, the degree of (in)visibility varies and depends in different ways on the sources outside and inside the cloaked region. These facts are somewhat surprising. The material in this section is taken from \cite{MinhLoc}. 
 In the last section, we discuss approximate cloaking for electromagnetic waves in the real time domain. We again implement the approach used for the acoustic waves. A new difficulty arises with this approach due to the fact that the multiplier technique, which played an important role in the analysis for the acoustic setting, does not quite suffice for the electromagnetic setting in the very high frequency regime. This point of difficulty is overcome by  a duality argument. The material in this section is taken from \cite{MinhLocT}. 

It is worth mentioning that there are other (than transformation) techniques which may be used to produce cloaking effects. Some schemes use negative index materials and complementary media, see, e.g.,   \cite{LaiChenZhangChanComplementary, Ng-Negative-Cloaking},  and some schemes obtain cloaking via localized resonance, see,  e.g., \cite{Ng-CALR-O} (and see also \cite{MiltonNicorovici, Ng-CALR}).   A brief survey on this topic can be found in \cite{Ng-Survey}.

\section{Approximate cloaking for acoustic waves in the time harmonic regime} \label{sec-2}

This section is devoted to approximate cloaking for acoustic waves in the time harmonic regime. The phenomena are thus  modelled by the  Helmholtz equation.  Our starting point is the regularization scheme introduced in \cite{Kohn}, in which one uses a transformation that blows up a small ball $B_\rho$ ($0< \rho < 1/2$) to the cloaked region $B_1$ in $\mR^d$ ($d=2, 3$). Here and in what follows $B_r$, $r > 0$,  denotes the ball centered at the origin of $\mR^d$ and of radius $r$. Our radial assumption on the geometry of the cloaked region is only for simplicity of notation. Consider the map  $F_{\rho}: \R^d \rightarrow \R^d$ defined by
\begin{align}\label{def-Frho} F_{\rho} (x)= \left\{ \begin{array}{cl} x &\text{ in } \R^d \setminus B_2, \\[6pt] 
\dsp  \left( \frac{2 - 2\rho}{2-\rho} + \frac{|x|}{2-\rho} \right)\frac{x}{|x|} &\text{ in } B_2 \setminus B_{\rho},\\[6pt]
\dsp \frac{x}{\rho} &\text{ in } B_{\rho}.\end{array} \right. \end{align}
The associated cloaking device  in the region $B_2 \setminus B_1$,  is characterized by the pair
\begin{equation}\label{TO--S}
\big( {F_{\rho}}_*I, {F_{\rho}}_*1 \big) \mbox{ in } B_2 \setminus B_1, 
\end{equation}
consisting of a matrix-valued function  and a scalar function. These are the material parameters in $B_2 \setminus B_1$, used in the corresponding  Helmholtz equation (see \eqref{eq-uc-1}).  Here and in what follows, we use the standard notation
\begin{equation}\label{pushforward}
F_*A(y) = \frac{\nabla F (x) A(x) \nabla F^T(x)}{ | \det \nabla F(x) |}, \quad F_*\Sigma(y) = \frac{\Sigma(x)}{ | \det \nabla F(x) |}, \quad x = F^{-1}(y),
\end{equation}
for the ``pushforward" of a symmetric, matrix-valued function $A$, and a scalar function $\Sigma$,  by the diffeomorphism $F$. $I$ denotes the identity matrix. 

We first consider the setting with a fixed lossy layer between 
the transformation based cloak and the cloaked region. We  assume that  the cloaked  region is $B_{1/2}$ to fix ideas, and we  assume that the lossy layer  in $B_1 \setminus B_{1/2}$ is characterized by the pair
\begin{equation}\label{F-lossy-layer}
\big(I, 1 + i / \omega \big), 
\end{equation} 
where $\omega > 0$ is the frequency and $i$ is the standard imaginary unit. We suppose that the medium outside $B_2$ (the cloaking device and the cloaked region) is homogeneous and the cloaked region is characterized by a pair $(a, \sigma)$ where 
$a$ is a real matrix-valued function and $\sigma$ is a real function, both defined in $B_{1/2}$. The medium in the whole space is thus given by  
\begin{equation}\label{med-thm-HF-1}
(A_c, \Sigma_c) = \left\{\begin{array}{cl} (I, 1) & \mbox{ in } \mR^d \setminus B_2, \\[6pt]
\big( {F_{\rho}}_*I, {F_{\rho}}_*1 \big)  & \text{ in } B_2 \setminus  B_1,\\[6pt]
I, 1 + i / \omega  & \text{ in } B_1 \setminus  B_{1/2},\\[6pt]
(a, \sigma)  &\text{ in } B_{1/2}.
\end{array} \right. 
\end{equation}
A related,  but slightly different setting was studied    in \cite{Kohn1}, where the authors  used  a $\rho$-dependent lossy layer. 
In that paper the lossy layer  in $B_1 \setminus B_{1/2}$ is characterized by the pair
\begin{equation}\label{F-lossy-layer}
\big(\rho^{d-2}, \rho^{d} + \rho^{d-2}i / \omega \big) = \big({F_{\rho}}_* I, {F_{\rho}}_*\left(1+i/\rho^2 \omega\right)   \big)~, 
\end{equation} 
and the cloaking problem is formulated in a bounded domain, but otherwise the medium is identical to (\ref{med-thm-HF-1}).
In all of the following it is assumed that $a$ is symmetric and uniformly elliptic, i.e.,   
\begin{equation}\label{elliptic}
\Lambda^{-1}|\xi|^2 \le \langle a(x)\xi, \xi \rangle \le \Lambda|\xi|^2\quad  \mbox{ for all } \, \xi \in \mR^d~,
\end{equation}
for almost every $x \in B_{1/2}$,
for some $\Lambda \ge 1$,  $\sigma$ is a positive function bounded above and below by positive constants.

Given a function $f \in L^2(\mR^d)$ and the frequency $\omega > 0$, the medium characterzied by $(A_c, \Sigma_c)$, produces the  unique outgoing  solution $u_c \in H^1_{\loc}(\mR^d)$ to the equation 
\begin{equation}\label{eq-uc-1}
\dive (A_c \nabla u_c) + \omega^2 \Sigma_c u_c = f \mbox{ in } \mR^d~.
\end{equation}
The homogeneneous medium characterized by $(I, 1)$, similarly produces the unique outgoing  solution $u \in H^1_{\loc}(\mR^d)$ to the equation 
\begin{equation}\label{eq-u}
\Delta u + \omega^2 u = f \mbox{ in } \mR^d~. 
\end{equation}
We refer to a solution $u \in H^1_{\loc}(\mR^d \setminus B_R)$  to the equation $\Delta u + \omega^2 u = 0 $ in $\mR^d \setminus B_R$  for some $R>0$, as outgoing, provided it satisfies
\begin{equation*}
\partial_r u - i \omega u = o( r^{\frac{1-d}{2}}) \mbox{ as } r = |x| \to + \infty~. 
\end{equation*}

The sense in which cloaking is achieved for a fixed frequency $\omega$ is described by the following result.
 
\begin{theorem} \label{thm-HF-1} \cite[Theorem 1.2]{Ng-1}, see also \cite[Theorem 3.1]{Kohn1}. Let $d=2, 3$, $0< \rho < 1/2$,  $R_0  > 2$,  and let $f \in L^2(\mR^d)$ with $\supp f \subset B_{R_0} \setminus B_2$. Assume that $u_c$ and $u$ are the unique outgoing solutions to \eqref{eq-uc-1} and \eqref{eq-u} respectively where $(A_c, \Sigma_c)$ is given in \eqref{med-thm-HF-1}. 
We then have 
\begin{equation*}
\| u_c - u\|_{H^1(B_R \setminus B_{2})} \le C e_d (\rho) \| f\|_{L^2(\mR^d)} ~\mbox{ for } R > 2~,  
\end{equation*}
for some positive constant $C = C_{R, \omega}$ independent of $f$, $\rho$, $\Lambda$, $a$, and $\sigma$. Here 
$$
e_d(\rho) = \rho \mbox{ if } d = 3 \mbox{ and } |\ln \rho|^{-1} \mbox{ if } d=2~. 
$$
\end{theorem}

As a consequence of  Theorem~\ref{thm-HF-1}, $\lim_{\rho \to 0} u_c = u$ in $\mR^d \setminus B_2$ for all $f$ with compact support outside $B_2$.  One  can therefore not detect the difference between $(A_c, \Sigma_c)$ and $(I, 1)$ as $\rho \rightarrow 0$ by observation of $u_c$ outside $B_2$: cloaking is achieved for observers outside $B_2$ in the limit as $\rho \rightarrow 0$.   

We now briefly describe the idea of the proof. As already mentioned, an essential ingredient of transformation based cloaking is the invariance of the Helmholtz equation under change of variables. This invariance may be stated as follows. 

\begin{lemma}\label{lem-TO-H}
Suppose $d \ge 2$. Let $A$ be a bounded  measurable, real matrix-valued function and $\Sigma$  be a bounded  measurable  complex function defined on $\mR^d$.  Let $F: \mR^d \mapsto \mR^d$ be bijective, with $F$ and $F^{-1}$ Lipschitz, and $\det \nabla F(x) > c>0$ for a.e. $x \in \mR^d$. Suppose $f \in L^2(\mR^d)$.  Then $u \in H^1_{\loc}(\mR^d)$ is a solution of
\begin{equation*}
\dive (A \nabla u) + \omega^2 \Sigma u = f \quad \mbox{ in } \mR^d
\end{equation*}
if and only if $v: = u \circ F^{-1} \in H^1_{\loc}(\mR^d)$ is a solution of
\begin{equation*}
\dive (F_*A \, \nabla v) + \omega^2 F_*\Sigma \,  v = F_*f \quad \mbox{ in } \mR^d.
\end{equation*}
\end{lemma}

\noindent Set 
\begin{equation*}
u_\rho = u_c \circ F_\rho \mbox{ in } \mR^d. 
\end{equation*}
Since $\supp f \cap B_2 = \emptyset $, it follows from Lemma~\ref{lem-TO-H} that  $u_\rho \in H^1_{\loc}(\mR^d)$ is the unique outgoing solution of 
\begin{equation}\label{eq-u-rho}
\dive(A_\rho \nabla u_\rho)  + \omega^2 \Sigma_\rho u_\rho = f \mbox{ in } \mR^d, 
\end{equation}
where 
\begin{equation*}
(A_\rho, \Sigma_\rho) = \left\{\begin{array}{cl} (I, 1) & \mbox{ in } \mR^d \setminus B_\rho, \\[6pt]
\big(\rho^{2-d}I, \rho^{-d}(1 + i / \omega)  \big) & \text{ in } B_\rho \setminus  B_{\rho/2},\\[6pt]
\big(\rho^{2-d}a(\cdot/ \rho), \rho^{-d}\sigma(\cdot/ \rho) \big)  &\text{ in } B_{\rho/2}.
\end{array} \right. 
\end{equation*}
Since $F_\rho(x) = x$ in $\mR^d \setminus B_2$, it is clear that 
\begin{equation*}
u_c - u = u_\rho - u \mbox{ in } \mR^d \setminus B_2. 
\end{equation*}
By comparing the coefficients of the equations satisfied by $u$ and $u_\rho$, one realizes that  the study of the cloaking effect can be reduced to the study of the effect of a small inhomogeneity.  
The effect of small Helmholtz inhomogeneities is well studied when the coefficients inside the small inclusions are fixed (or have a finite range), see,  e.g., \cite{FV, VogeliusVolkov}. 
Nevertheless, the situation in the cloaking context is non-standard since the coefficients inside the small inclusion blow up as the diameter goes to 0. To deal with this, we split $u_\rho - u$ into two parts $w_{1, \rho}$ and $w_{2, \rho}$ where 
\begin{equation*}
w_{1, \rho} = u_{1, \rho} - u  \mbox{ in } \mR^d \quad \mbox{ and } \quad w_{2, \rho} = u_{\rho} - u_{1, \rho}  \mbox{ in } \mR^d . 
\end{equation*}
Here $u_{1, \rho} \in H^1_{\loc}(\mR^d)$ is the unique outgoing  solution of the system 
\begin{equation}\label{def-u1rho}
\left\{\begin{array}{cl}
\Delta u_{1, \rho} + \omega^2 u_{1, \rho} = f  & \mbox{ in } \mR^d \setminus B_\rho, \\[6pt]
u_{1, \rho}  = 0 & \mbox{ in } B_\rho. 
\end{array}\right. 
\end{equation}
In other words, $u_{1, \rho}$ satisfies the Dirichlet problem with zero boundary in the exterior domain $\mR^d \setminus B_\rho$ and $u_{1, \rho} = 0$ in $B_{1, \rho}$.  This way of  splitting  is inspired by \cite{Ng-Vogelius-1} in which uniform estimates for the scattering effects of small (conductivity) inhomogeneities were studied.  The functions $w_{1,\rho}$ and $w_{2,\rho}$ satisfy
$$
\left\{\begin{array}{cl}
\Delta w_{1, \rho} + \omega^2 w_{1, \rho} = 0  & \mbox{ in } \mR^d \setminus B_\rho, \\[6pt]
w_{1, \rho}  = -u & \mbox{ on } \partial B_\rho~, 
\end{array}\right. 
$$
and
\begin{equation*}
\left\{\begin{array}{ll}
\Delta w_{2,\rho} + \omega^2 w_{2,\rho} = 0 & \mbox{in } \mR^d \setminus B_\rho, \\[6pt]
\dive (A_\rho \nabla w_{2,\rho}) + \omega^2 \Sigma_\rho w_{2,\rho} = 0 & \mbox{in } B_\rho, \\[6pt]
\dsp \frac{\partial w_{2,\rho}}{\partial \nu}\Big|_{\mathrm{ext}} - \frac{1}{\rho^{d-2}} \frac{\partial w_{2,\rho}}{\partial \nu}\Big|_{\mathrm{int}} =  -\frac{\partial u_{1,\rho}}{\partial \nu}\Big|_{\mathrm{ext}}& \mbox{on } \partial B_\rho~, 
\end{array}\right.
\end{equation*}
respectively. By  a change of variables $x \mapsto x/ \rho$ and subsequent analysis, one obtains the following estimates for $w_{1, \rho}$ and $w_{2, \rho}$
\begin{equation*}
\| w_{1, \rho}\|_{L^2(B_R \setminus B_2)} \le C e_d(\rho) \|f \|_{L^2}
\end{equation*}
and 
\begin{equation*}
\| w_{2, \rho}\|_{L^2(B_R \setminus B_2)} \le C \rho^{d-1} \|f \|_{L^2}~. 
\end{equation*}
These together yield the desired estimate for $u_\rho -u$ ($=u_c-u$ outside $B_2$). Just a few comments on the above estimates for $w_{1,\rho}$ and $w_{2,\rho}$:
the proof of the first estimate (on $w_{1,\rho}$) is quite standard in three dimensions, thanks to the ``boundedness" of the fundamental solution of the Helmholtz equation with respect to frequency (as it approaches zero). The proof in two dimensions is more involved  due to the fact that  the fundamental solution blows up as the frequency goes to 0.  To handle this situation, we ``decompose" the solution into two parts: the first one accounting for the behavior of the fundamental solution and then a remainder  which is easier to  control. 
The proof of the second estimate (on $w_{2,\rho}$) is based on a variational argument and proceeds by contradiction.
The presence of the lossy layer plays a decisive role in the proof of this second estimate.  The lossy layer allows one to obtain a control on  $\| w_\rho\|_{L^2(B_\rho \setminus B_{\rho/2})}$ essentially by multiplying  
the equation by $\bar w_{2,\rho}$ (the conjugate of $w_{2,\rho}$), integrating on $B_R$, letting $R$ go to infinity, and considering the imaginary part as usual. We refer the reader to \cite{Ng-1} for more details.

In \cite{Ng-Vogelius-2}, we studied how the degree of visibility depends on frequency.  Instead of a fixed lossy layer we employed a lossy layer depending on $\rho$: 
\begin{equation}\label{V-lossy-layer}
{F_{\rho}}_*I, {F_{\rho}}_* \left( 1 + \frac{i }{\omega \rho \lambda}\right) \mbox{ in } B_1 \setminus B_{1/2}, 
\end{equation}
where $0< \lambda< 1$ is a fixed parameter.  Otherwise the medium is as in \eqref{F-lossy-layer}.

The estimates obtained in the low frequency regime rely on the behavior of the fundamental solution of the Helmholtz equation in this regime. These estimates  are weaker than the ones in the high frequency regime. The results in \cite{Ng-Vogelius-2} are optimal and compatible with the ones in Theorem~\ref{thm-HF-1}.   The strategy of the analysis is similar to the one of Theorem~\ref{thm-HF-1} by deriving estimates on the effect of  small inhomogeneities, but now the dependence on frequency is explicit. An important ingredient in the analysis in the  high frequency regime  is the multiplier technique. This technique has  its roots in the work of Morawetz and Ludwig  \cite{MorawetzLudwig} (see also the work of Rellich \cite{Rellich},  and Perthame and Vega \cite{PerthameVega}). The lossy layer plays an important role in the analysis: without it, our estimates might not hold, due to the trapping phenomenon associated with the Helmholtz equation in the high frequency regime.  The details of the analysis is outside the scope of this review and can be found in \cite{Ng-Vogelius-2}. 

We next turn to the situation where no lossy layer is employed. We assume here that  the cloaked region is $B_1$. The cloaked object is characterized by the pair  $(a, \sigma)$, which is  assumed to occupy $B_1$ and satisfy the same assumptions as before.  With the cloaking device and the cloaked object, the medium is 
\begin{equation}\label{med-thm-HF-2}
(A_c, \Sigma_c) = \left\{\begin{array}{cl} (I, 1) & \mbox{ in } \mR^d \setminus B_2, \\[6pt]
\big( {F_{\rho}}_*I, {F_{\rho}}_*1 \big)  & \text{ in } B_2 \setminus  B_1,\\[6pt]
(a, \sigma)  &\text{ in } B_{1}.
\end{array} \right. 
\end{equation}
We  consider now also the case where $f$ does {\it not necessarily} vanish in the cloaked region $B_1$, however, $f$ is still assumed to be 0 in the ``cloaking device" region $B_2 \setminus B_1$.  An important concept that now appears is the concept of resonance. In three dimensions, this concept  is related to the Helmholtz equation inside $B_1$ together with the zero Neumann boundary condition:
\begin{definition} \cite{Ng-2} 
Let $d=3$ and define 
\begin{equation}\label{defMB1}
{\mathcal M}: = \Big\{\psi \in H^1(B_1): \; \dive (a \nabla \psi) + \omega^2 \sigma \psi = 0 \mbox{ in } B_1 \mbox{ and } a \nabla \psi \cdot \nu = 0 \mbox{ on } \partial B_1 \Big\}.
\end{equation}
The system is called non-resonant if  ${\mathcal M} = \{0 \}$, otherwise, the system is called resonant. 
\end{definition}
We have

\begin{theorem} \cite[Theorem 1.4 and Proposition 1.11]{Ng-2} 
 Let $d=3$, $\omega > 0$, and $0< \rho < 1/2$, $R_0>2$, and let $f \in L^2(\mR^3)$ with $\supp f \subset (B_{R_0} \setminus B_2)\cup \overline{B_1}$. 
Assume that $u_c$ and $u$ are the unique outgoing solutions of \eqref{eq-uc-1} and \eqref{eq-u} respectively, in the latter case with $f$ replaced by $f\mathds{1}_{\{|x|\ge 2\}}$. The coefficient pair $(A_c, \Sigma_c)$ is given by \eqref{med-thm-HF-2}. The following statements hold
\begin{enumerate}
\item Assume  that $\dsp \int_{B_1} f \bar \be = 0$ for all $\be \in {\mathcal M}$. Then for all $K \subset \subset \mR^3 \setminus \overline B_1$,
\begin{equation*}
\| u_c - u \circ F_0^{-1} \|_{H^1(K)}  \le C \rho \| f \|_{L^2},
\end{equation*}
for some positive constant $C$ depending on $\omega$, $K$, $a$, and $\sigma$, but independent of $\rho$, $f$. Moreover,   
$
u_c \mbox{ converges weakly in } H^1(B_1)
$
as $\rho \to 0$. 

\item Assume that ${\mathcal M} \not =\{ 0 \}$ and $f = \be \mathds{1}_{B_1}$ with $\be \in {\mathcal M} \setminus \{0 \}$. Then 
\begin{equation*}
\liminf_{\rho \to 0}  \rho \| u_c\|_{H^1(B_1)} > 0.
\end{equation*}
Assume in addition that $\be$ is radial, $a$ and $\sigma$ are isotropic and homogeneous in $B_1$,  i.e.,   $a = \lambda_1 I$ and $\sigma= \lambda_2$ for some positive constants $\lambda_1$ and $\lambda_2$. Then
\begin{equation*}
\liminf_{\rho \to 0} \| u_c\|_{L^2(B_4 \setminus B_2)} > 0.
\end{equation*}
\end{enumerate}
\end{theorem}

Here and in what follows, $\mathds{1}_D$ denotes the characteristic function of a subset $D$ of $\mR^d$. 

\medskip 

The limit of $u_c$ in $B_1$ is well-determined in the first assertion. In the non-resonant case ${\mathcal M} = \{0\}$, the limit is just the solution to the corresponding Neumann problem, which is unique. In the resonant case ${\mathcal M} \not = \{0\}$, the limit is uniquely determined by a (non-local) transmission problem which involves the values of $u(0)$ and the corresponding Neumann problem, see \cite[Definition 1.5]{Ng-2} for the details. 
The cloaking estimates  are compatible with Theorem~\ref{thm-HF-1}, but note that the constant $C$ depends on $a$ and $\sigma$.  The non-cloaking and infinite energy results in the resonant case are new and quite surprising. In the literature, only the concept of finite energy had been  considered  for the singular scheme, see,  e.g., \cite{GKLU07-1}.


In the two dimensional case, the definition of resonance is quite different and depends on  a non-local transmission problem. 
The different definitions of resonance in two and three dimensions arise from the difference in  $(A_\rho, \Sigma_\rho)$ in two and three dimensions. 
We have 

\begin{definition} \cite[Definition 1.7]{Ng-2} 
Let $d=2$.  Define 
$$
{\mathcal N} = \{ \psi \in W^{1}(\mR^2) \mbox{ such that \eqref{uniqueness2d} holds}\}, 
$$
where 
\begin{equation}\label{uniqueness2d}
\left\{\begin{array}{ll}
\Delta \psi= 0 & \mbox{in } \mR^2 \setminus B_1, \\[6pt]
\dive (a \nabla \psi) + \omega^2 \sigma \psi = 0 & \mbox{in } B_1, \\[6pt]
\dsp \frac{\partial \psi}{\partial \nu}\Big|_{\mathrm{ext}}  =   a \nabla \psi \cdot \nu  \Big|_{\mathrm{int}} & \mbox{on } \partial B_1.
\end{array}\right.
\end{equation}
The system is  called non-resonant if ${\mathcal N} = \{0 \}$,  otherwise, the system is called resonant. 
\end{definition}

Here  
\begin{equation*}
W^{1}(\mR^2) = \Big\{ \psi \in L^1_{loc}(\mR^2): \frac{\psi(x)}{\ln(2 + |x|) \sqrt{1 + |x|^2}}  \in L^2(\mR^2) \mbox{ and } \nabla \psi \in L^2(\mR^2) \Big\}.
\end{equation*}

We have 
\begin{theorem}   \cite[Theorem 1.8 and Proposition 1.12]{Ng-2} 
Let $d=2$, $\omega > 0$, $0 < \rho < 1/2$, $R_0>2$, and let $f \in L^2(\mR^3)$ with $\supp f \subset (B_{R_0} \setminus B_2)\cup \overline{B_1}$. Assume that $u_c$ and $u$ are the unique outgoing solutions of \eqref{eq-uc-1} and \eqref{eq-u} respectively, in the latter case with $f$ replaced by $f  \mathds{1}_{\{ |x|\ge 2\}} $. The coefficient pair  $(A_c, \Sigma_c)$ is given by \eqref{med-thm-HF-2}. The following statements hold 

\begin{enumerate}
\item Assume that the system is non-resonant. Then for all $K \subset \subset \mR^2 \setminus \overline B_1$,
\begin{equation*}
\| u_c - u \circ F_0^{-1} \|_{H^1(K)}  \le \frac{C}{|\ln \rho|} \| f \|_{L^2},
\end{equation*}
for some positive constant $C$ depending on $\omega$,  $K$,  $a$, and $\sigma$, but independent of $\rho$ and $f$. Moreover,  $u_c$ weakly converges  in $H^1(B_1)$ as $\rho \to 0$. 

\item Assume that the system is resonant. Fix  $\be \in {\mathcal N} \setminus \{0 \}$ and set  $f = \mathds{1}_{B_1} \be$.  Then 
$$
\lim_{\rho \to 0} \| u_c\|_{H^1(B_1)} = + \infty. 
$$
Assume in addition that $\be$ is radial, $a$ and $\sigma$ are isotropic and homogeneous in $B_1$.  Then
$$
\liminf_{\rho \to 0 } \| u_c\|_{L^2(B_4 \setminus B_2)} = + \infty. 
$$

\end{enumerate}
\end{theorem}

The facts that resonance is defined in terms of a non-local problem, that cloaking might not be achieved in the resonant case and that the energy of the fields inside the cloaked region can blow up are  new and somewhat surprising.  The limit of $u_c$ in $B_1$ is well-determined in the non-resonant case, see \cite[Definition 1.7]{Ng-2}. 

We end this section with the topic of enhancement of cloaking. Its motivation comes from the desire to reduce the complexity of the materials needed to achieve a prescribed level of cloaking.  In \cite{GV}, with Griesmaier, the second author introduce, for the electrostatic case, a regularized, approximate cloaking by mapping scheme and analyse the problem of optimal choice among radial maps. Two different optimality criteria are investigated: minimal maximal anisotropy and minimal mean anisotropy of the conductivity distribution. Using both criteria, it is shown that it is possible to achieve significantly lower anisotropy (for a prescribed level of invisibility) or significantly lower visibility (for a prescribed level of anisotropy). For example, in two dimensions one may achieve exponentially small visibility with a cloak, that in terms of anisotropy (and lowest and highest conductivity) is no worse than the traditional affine map cloak, which only yields quadratically small visibility, see, e.g., \cite{Kohn, Ng-Vogelius-1} and \cite{HV}. A discussion of the  sense in which  radial maps are "best possible" if the cloak occupies $B_2 \setminus B_1$ is found in \cite{Cap-Vog}.

\section{Approximate cloaking for acoustic waves in the time regime} \label{sec-3}

This section is devoted to approximate cloaking in the time domain. The material  presented here is from \cite{Ng-Vogelius-3} and \cite{ Ng-Vogelius-4} which are,  to the best of our knowledge, the first mathematical works on approximate cloaking for the full wave equation. 

We first discuss the cloaking scheme used in \cite{Ng-Vogelius-3}. The cloaking device itself has two layers. In frequency domain this cloaking device is  given by \eqref{TO--S} and \eqref{V-lossy-layer} with $\lambda = \rho^{1+\gamma}$. More precisely, we have
\begin{equation*}
\big(A_c, \Sigma_{1, c}, \Sigma_{2, c} \big) = \left\{\begin{array}{cl} \big(I, 1, 0 \big) & \mbox{ in } \mR^d \setminus B_2, \\[6pt]
\big( {F_{\rho}}_*I, {F_\rho}_*1, 0 \big) & \mbox{ in } B_2 \setminus B_1, \\[6pt]
\big( {F_{\rho}}_*I, {F_\rho}_*1, {F_{\rho}}_*(1/ \rho^{2+\gamma})  \big)& \mbox{ in } B_1 \setminus B_{1/2}, \\[6pt]
\big(a, \sigma, 0 \big) & \mbox{ in } B_{1/2},  
\end{array}\right. 
\end{equation*}
where $\Sigma_{1, c}$ denotes the real part of $\Sigma_c$ and $\Sigma_{2, c}$ denotes the imaginary part of $\omega \Sigma_c$. In time domain the cloaking scheme thus gives rise to the system 
\begin{equation}\label{w-uc}
\left\{\begin{array}{cl}
\Sigma_{1, c}\partial_{tt}^2 u_c- \dive(A_c \nabla  u_c) + \Sigma_{2, c} \partial_t u_c = f  & \mbox{ in } (0, \infty) \times \mR^d, \\[6pt]
u_c(t =0, \cdot)  = \partial_t u_c(t =0, \cdot) = 0 \mbox{ in } \mR^d. 
\end{array}\right. 
\end{equation}
The homogeneous medium in time domain corresponds to the system
\begin{equation}\label{w-u}
\left\{\begin{array}{cl}
\partial_{tt}^2  u - \Delta u = f & \mbox{ in } (0, \infty) \times \mR^d, \\[6pt]
u(t =0, \cdot)  = \partial_t u(t =0, \cdot)  = 0 &  \mbox{ in } \mR^d.
\end{array}\right. 
\end{equation}
We   consider here the situation where the initial time is 0 and the initial conditions are zero. 
Given $f \in L^2\big((0, \infty); L^2(\mR^d) \big)$, both systems (\ref{w-uc}) and (\ref{w-u}) have unique (weak) solutions $u_c, u \in  L^2_{\loc}\big([0, \infty); H^1(\mR^d) \big)$ with $\partial_t u_c, \partial_t u \in L^\infty_{loc}\big([0, \infty); L^2(\mR^d) \big)$.

\medskip 
The approximate cloaking property of this scheme is described by

\begin{theorem}   \cite[Theorems 1 and 2]{Ng-Vogelius-3} \label{thm2-HF}  Let $0 < \rho < 1/ 2$, $T_0>0$,  and $R_0>2$, and let $f \in   L^2\big((0, \infty); L^2(\mR^d) \big)$. Assume that 
$ \supp f \subset (0, T_0) \times (B_{R_0} \setminus B_2)$ and let $u_c$ and $u$ be the unique solution to \eqref{w-uc} and \eqref{w-u} respectively . We have, for $R>2$, 
\begin{equation}\label{est-thm2-HF}
\sup_{t > 0} \| u_c(t, \cdot) - u(t, \cdot)\|_{L^2(B_R \setminus B_2)} \le C e_d(\rho) \| f\|_{C^m \big((0,T_0) \times \mR^d \big)}, 
\end{equation}
for some $m>0$. Here $C=C_R$ is a positive constant depending on the ellipticity of $a$, the upper and lower bounds of $\sigma$, $R_0$, and $T_0$ but independent of $\rho$, $f$, $a$, and $\sigma$. 
\end{theorem}

To obtain these estimates of the degree of near-invisibility for the full wave equation, we proceed as follows. We first transform the wave equation into a family of Helmholtz equations by taking the Fourier transform with respect to time. After obtaining the appropriate degree of near-invisibility estimates for the Helmholtz equations, where the dependence on frequency is explicit, we simply invert the Fourier transform. For the high frequency regime we can directly use the estimates on the degree of near-invisibility established in \cite{Ng-Vogelius-2}. However, for the low frequency regime, we  establish new estimates which improve on the dependence on the frequency compared to the ones in \cite{Ng-Vogelius-2}. It is at this place that we use the (additional) finite range assumption on $(a, \sigma)$. The proof in places uses the theory of $H$-convergence (see, e.g., \cite{MuratTartar}).  Our estimates in the frequency domain blow up as the frequency goes to 0 but  they blow up in an integrable way as established in \cite[Section 2.2]{Ng-Vogelius-3}.  This is sufficient to yield the estimates in time domain. 

As an important technical point we need to establish that the Fourier transforms of solutions to the wave equation (with respect to time) satisfy the outgoing radiation conditions. This fact is of independent interest, and  very useful in the study of various problems in time domain, since it allows one directly to make use of  knowledge from the frequency domain, see, e.g.,   \cite{Ng-Vogelius-4} and \cite{MinhLinh}.

In \cite{Ng-Vogelius-4}, we investigated approximate cloaking for a wave equation in which the transformation based cloak obeys the Drude-Lorentz model.
More precisely, instead of \eqref{TO--S}, we used the more physically realistic model
\begin{equation*}\label{TO-S-L}
{F_\rho}_*I, {F_\rho}_* 1 + \sigma_{1, c} (\omega) \mbox{ in } B_2 \setminus B_1, 
\end{equation*}
where 
\begin{equation*}
\sigma_{1, c}(\omega) = \frac{\sigma_N}{\omega_c^2 - \omega^2 - i \sigma_D \omega}
\end{equation*}
Here $\sigma_N$ and $\sigma_D$ are material constants which can in principle depend on the space variable $x$, and $\omega_c > 0$ is the so-called resonant frequency of the Drude-Lorentz model.  We still employed a fixed lossy layer in $B_1 \setminus B_{1/2}$.  Under the assumption that the resonant frequency $\omega_c$ satisfies
$c_*\rho^{-d/2} \le \omega_c\le C_* \rho^{-K}~,$
we showed that cloaking is achieved. However, instead of 
\eqref{est-thm2-HF}, one has 
\begin{equation}\label{est-thm2-HF-v}
\sup_{t \in (0, T)} \| u_c(t, \cdot) - u(t, \cdot)\|_{L^2(B_R \setminus B_2)} \le C e_d(\rho) \| f\|_{C^m}, 
\end{equation}
for some $m>0$. Here $C = C_{R, T}$ is a positive constant independent of $a$, $\sigma$, $f$, and $\rho$.  This result
is  in a slightly different spirit than the one in Theorem~\ref{thm2-HF}, since the constant $C$ here is independent of $a$ and $\sigma$.  This constant depends on $T$, though. 

The approach in \cite{Ng-Vogelius-4} borrows several ideas from our previous work  in \cite{Ng-Vogelius-3}. We transform the wave equation into a family of Helmholtz equations by taking the Fourier transform with respect to time. Having established  appropriate near-invisibility estimates for the Helmholtz equations, with explicit frequency dependence, we then essentially invert the Fourier transform.  
An  important aspect of the analysis is that one needs to deal with a wave equation which is non-local in time, due to $\sigma_{1, c}$'s dependence on frequency.  A key fact used  in our analysis is the causality of the Drude-Lorentz model. This fact has later been used to establish the well-posedness  of the Maxwell equations in the time domain  for very general  dispersive materials \cite{Ng-Vinoles}. 
We also use a helpful idea from \cite{MinhLinh} to simplify the analysis (and avoid the problems of low frequency blow-up). 
The idea is to initially control $\partial_t u_c - \partial_t u$, using the the frequency domain estimates in \cite{Ng-Vogelius-2}, and then turn this into control of   $u_c -  u$, on a finite interval of time,  by integration. It is the use of this idea which causes the constant $C$ to depend on $T$ (but stay independent of $a$ and $\sigma$). 

\section{Approximate cloaking for electromagnetic waves in the time harmonic regime} \label{sec-4}

This section is devoted to approximate cloaking for Maxwell equations in the time harmonic regime using  schemes in the spirit of \cite{Kohn}. Other researchers have worked  in this direction. In \cite{Ammari13}, Ammari et al. investigated cloaking using appropriate additional  layers between the transformation cloak and the cloaked region in order to  enhance the cloaking effect. Their technique is based on knowledge of the polarization tensors and requires the materials inside the cloaked region to be constant and isotropic. In \cite{Bao}, Bao, Liu, and Zou  studied approximate cloaking using  a lossy layer in the spirit of what was described  in the previous sections for the Helmholtz equation. Their approach is based on standard estimates in the lossy layer and does not provide an optimal estimate of the degree of visibility.   
Using separation of variables, Lassas and Zhou \cite{Lassas}, studied the transformation cloak in a symmetric setting and examined the limit of solutions to the approximate cloaking problem near the cloak interface, in the  non-resonant case (see Definition~\ref{def-R}). The study of finite energy solutions for the singular scheme can be found in \cite{GKLU07-1, Weder1, Weder2}.

We now desribe some results  due to Tran and the first author, \cite{MinhLoc}. We consider the situation where 
the cloaking device {\it only} consists of a layer constructed by the mapping technique. Due to the fact that no lossy  layer is present, resonance might appear and therefore the analysis is quite  subtle and the phenomena are complex.  
We investigate both the non-resonant  and the resonant case (Definition~\ref{def-R}). The results reveal several new phenomena in comparison with the works mentioned above, and in comparison with the acoustic setting as well. 

To fix notation, we suppose the cloak occupies
the  region $B_2 \setminus B_1$ and the cloaked region is the unit ball $B_1$  inside which the permittivity and the permeability are given by two $3\times 3$ matrix-valued functions $\eps_O$ and $ \mu_O$ respectively.  
From now on, we assume that 
\begin{equation}\label{symetric}
\eps_O, \,  \mu_O \mbox{ are real, symmetric},  
\end{equation}
and uniformly elliptic in $B_1$. 
We also assume that $\eps_O, \mu_O$ are piecewise $C^1$ in order to ensure the well-posedness of Maxwell's equations in the frequency domain via the unique continuation principle (see, e.g., \cite{Protter, Nguyen, BallCapdeboscq}). In the spirit of  the scheme in \cite{Kohn}, the pair of  permittivity and permeability of the cloaking region is given by
\begin{equation}\label{Frho-M}
 ({F_{\rho}}_*I, {F_{\rho}}_*I ) \mbox{ in } B_{2} \setminus B_1, 
\end{equation}
where $F_{\rho}$ is given by \eqref{def-Frho}. We assume the medium is homogeneous outside the cloak and the cloaked region. In the presence of the cloaked object and the cloaking device, the medium in the whole space $\R^3$ is then given by  $(\eps_c, \mu_c)$  defined as follows
\begin{equation}
\label{medium:cloak}
(\eps_c, \mu_c) = \left\{\begin{array}{cl} (I, I) & \mbox{ in } \mR^3 \setminus B_2, \\[6pt]
\big( {F_{\rho}}_*I, {F_{\rho}}_*I \big)  & \text{ in } B_2 \setminus  B_1,\\[6pt]
(\eps_O, \mu_O)  &\text{ in } B_{1}.
\end{array} \right. 
\end{equation}

In the medium characterized by $(\eps_c, \sigma_c)$, at the  frequency $\omega>0$, the electromagnetic field generated by a current $J \in [L^2(\mR^3)]^3$ is the unique (Silver-M\"uller) radiating solution $(E_c, H_c)\in [H_{\loc}(\curl, \R^3)]^2$ to the system
\begin{equation}
\label{equ:cloak}
\begin{cases}
\nabla \times E_c = i\omega \mu_c H_c  &\text{ in } \mathbb{R}^3,\\[6pt]
\nabla \times H_c = -i\omega \eps_c E_c + J &\text{ in } \mathbb{R}^3.
\end{cases}\end{equation}
For an open subset $U$ of $\mR^3$, we use the standard notations
\[
H(\curl, U) := \Big\{\phi \in [L^2(U)]^3:  \;  \nabla \times \phi \in [L^2(U)]^3 \Big\}
\]
and 
\[
H_{\loc}(\curl, U) := \Big\{\phi \in [L_{\loc}^2(U)]^3: \;  \nabla \times \phi \in [L_{\loc}^2(U)]^3 \Big\}.
\]
Recall that a solution $(E, H) \in [H_{\loc}(\curl, \R^3\setminus B_R)]^2$ to the Maxwell equations 
\[
\begin{cases}
\nabla \times E = i \omega H  &\text{ in } \mathbb{R}^3\setminus B_R,\\[6pt]
\nabla \times H = -i \omega E  &\text{ in } \mathbb{R}^3\setminus B_R
\end{cases}
\]
(for some $R> 0$) is called radiating if it satisfies the (Silver-Muller) radiation conditions
\begin{equation}\label{SM-condition}
H \times x - |x| E = O(1/|x|) \quad   \mbox{ and } \quad  E\times x + |x| H = O(1/|x|) \mbox{ as } |x| \to + \infty~. 
\end{equation}
 We assume that $J = 0$ in $B_{2} \setminus B_1$ (no source in the cloaking device).  Denote by  $J_{\ext}$ and $J_{\inte}$ the restriction of $J$ to $\mR^3 \setminus B_2$ and $B_1$, respectively. 
In the homogeneous medium (in the absence of the cloaking device and the cloaked object), the electromagnetic field  generated by $J_{\ext}$ is the unique (Silver-M\"uller) radiating solution $(E, H)\in [H_{\loc}(\curl, \R^3)]^2$ to the system 
\begin{align}
\label{equ:homo}
\begin{cases}
\nabla \times E = i\omega H  &\text{ in } \mathbb{R}^3,\\[6pt]
\nabla \times H = -i\omega E + J_{\ext}  &\text{ in } \mathbb{R}^3.
\end{cases}
\end{align}

We next introduce the notion of resonance.  

\begin{definition} \label{def-R}   Set 
 \[{\mathcal N} =  \Big\{ (\bE, \bH) \in [H(\curl, B_1)]^2:  (\bE, \bH) \text{ satisfies the system } (\ref{def:resonance}) \Big\}, \]
where
\begin{equation}
\label{def:resonance}
\begin{cases} \nabla\times \bE = i\omega \mu_O \bH & \text{ in } B_1,\\[6pt]
 \nabla\times \bH = -i\omega \eps_O \bE &\text{ in } B_1, \\[6pt] 
 \nabla\times \bE \ \cdot \ \nu = \nabla\times \bH \ \cdot \ \nu = 0 &\text{ on } \partial B_1.
 \end{cases}
\end{equation} 
The cloaking system \eqref{medium:cloak} is called non-resonant if ${\mathcal N} = \{(0, 0) \}$. Otherwise, the  cloaking system \eqref{medium:cloak} is called resonant. 
\end{definition}

Denote 
\begin{equation}\label{def-F*E}
   	F*E: = (\nabla F^{-T}E)\circ F^{-1}
\end{equation}
for a vector field $E$ and a diffeomorphism  $F$. 
Our main result in the non-resonant case is the following theorem. 
  
   \begin{theorem}
   \label{thm1}
   Let $0 < \rho < 1/2$,  $R_0> 2$, and  let $J \in L^2(\mR^3)$ be  such that $\operatorname{supp} J_{\ext} \subset \subset B_{R_0}\setminus B_2$.  Assume that system (\ref{medium:cloak}) is non-resonant. We have, for all $K\subset \subset \R^3\setminus \bar{B_1}$, 
 \begin{equation}\label{thm1-est1}
 \|(F_\rho^{-1}*E_c, F_\rho^{-1}*H_c) - (E, H)\|_{H(\curl, K)}\leq C\left(\rho^3 \|J_{\ext}\|_{L^2(B_{R_0}\setminus B_2)}+\rho^2\|J_{\inte}\|_{L^2(B_1)}\right),
\end{equation}
   for some positive constant $C$ depending only on $R_{0}, \omega, K,  \eps_O$, and $\mu_O$. Moreover, $(E_c, H_c)$ converges in $ [H(\curl, B_1)]^2$,  as $\rho \to 0$, to the unique solution $(\bE, \bH) \in [H(\curl, B_1)]^2$ of 
   \begin{equation*}
\begin{cases} \nabla\times \bE = i\omega \mu_O \bH & \text{ in } B_1,\\[6pt]
 \nabla\times \bH = -i\omega \eps_O \bE + J_{int} &\text{ in } B_1, \\[6pt] 
 \nabla\times \bE \ \cdot \ \nu = \nabla\times \bH \ \cdot \ \nu = 0 &\text{ on } \partial B_1.
 \end{cases}
\end{equation*}
   \end{theorem}

The novelty of Theorem~\ref{thm1}, in comparison with the works mentioned above, is as follows. (1) No lossy layer is present, and (2) the result holds for a general class of pairs $(\eps_O, \mu_O)$.  Applying the technique used to prove Theorem~\ref{thm1} to the case where a fixed lossy layer is present, one improves a result from \cite{Bao}, where the authors obtained an estimate for the far-field visibility in terms of the sum of $\rho^2$ times the magnitude of the incident field and $\rho$ times the norm of  the source in the cloaked region. In contrast to \cite{Ammari13, Bao},  Theorem~\ref{thm1} provides an estimate of visibility up to the cloaked region, and the behavior of the electromagnetic fields are established inside the cloaked region. Some comments are in order concerning the boundary conditions of the limit  in Theorem~\ref{thm1} in comparison with the works in \cite{GKLU07-1, Weder2} where the authors considered finite energy solutions for the singular scheme.   In \cite{Weder2}, the conditions 
$$
 \nabla\times E_c \cdot \nu |_{int} = \nabla\times H_c \cdot \nu |_{int}= 0
$$
are also imposed on the boundary of the cloaked region. However, this is different from \cite{GKLU07-1} (see also \cite[page 459]{Lassas}), where the following boundary  conditions are imposed
$$
E_c \times  \nu|_{int} = H_c \times  \nu|_{int} = 0. 
$$

\medskip 

We next consider the resonant case.

   \begin{theorem}
   \label{thm1.1} 
  Let $0< \rho < 1/ 2$,  $R_0> 2$, and let $J \in [L^2(\mR^3)]^3$ be  such that $\operatorname{supp} J_{\ext} \subset \subset B_{R_0}\setminus B_2$.  Assume that  system (\ref{medium:cloak}) is resonant. We have 
  \begin{enumerate}
  \item If  the following compatibility condition {\bf holds}: 
    \begin{equation}
   \label{passive}
   \int_{B_1}J_{\inte}\cdot\bar{\bE}\,dx = 0 \quad \mbox{ for all } (\bE, \bH) \in {\mathcal N}. 
   \end{equation}
then, for all $K\subset \subset \R^3\setminus \bar{B_1}$, 
 \begin{equation}\label{thm1.1-est1}
 \|(F_\rho^{-1}*E_c, F_\rho^{-1}*H_c) - (E, H)\|_{H(\curl, K)} \leq C \Big(\rho^3 \|J_{\ext}\|_{L^2(B_{R_0}\setminus B_2)}+\rho^2\|J_{\inte}\|_{L^2(B_1)} \Big), 
 \end{equation}
   for some positive constant $C$ depending only on $R_{0}, \omega, K, \mu_O$, and $\eps_O$. Moreover, $(E_c, H_c)$ converges in $ [H(\curl, B_1)]^2$ as $\rho \to 0$. 
   
   \item If the compatibility condition does {\bf not} hold, i.e.,      \begin{equation}
   \label{N-passive}
   \int_{B_1}J_{\inte}\cdot\bar{\bE}\,dx \neq  0 \quad \mbox{ for some } (\bE, \bH) \in {\mathcal N}, 
   \end{equation}
then, for all $K\subset \subset \R^3\setminus \bar{B_1}$,
\begin{equation}\label{thm1.2-est1}
\| (F_\rho^{-1}*E_c, F_\rho^{-1}*H_c) - (E, H)\|_{H(\curl, K)} \leq C \Big(\rho^3 \|J_{\ext}\|_{L^2(B_{R_0}\setminus B_2)}+\rho\|J_{\inte}\|_{L^2(B_1)} \Big)
\end{equation}
and
   \begin{equation} \label{explosion}\liminf \limits_{\rho\rightarrow 0}  \rho\| \big(E_c, H_c \big)\|_{L^2(B_1)} > 0.
   \end{equation}
   \end{enumerate}
  \end{theorem}
  
The limit of $(E_c, H_c)$ in the first part of Theorem~\ref{thm1.1} is also uniquely  determined. It solves a nonlocal problem, as in the case of the Helmholtz equation, see \cite{MinhLoc} for more details. 
\medskip

Some comments about Theorem~\ref{thm1.1} are in order.  Theorem~\ref{thm1.1} implies in particular that cloaking is achieved even in the resonant case. Moreover, without source in the cloaked region, one can achieve the same degree of (in)visibility as in the non-resonant case considered in Theorem~\ref{thm1}. In general, the degree of visibility varies and depends on the compatibility of the source inside the cloaked region. More precisely,   the rate of the convergence of $(E_c, H_c) - (E, H)$ outside $\bar B_1$ in the compatible  case is of the order $\rho^2$ which is better than in the incompatible  case, where an estimate of the order $\rho$ is obtained. These rates of convergence are optimal, as  discussed in \cite{MinhLoc}. By  \eqref{explosion}, in the incompatible case the energy inside the cloaked region blows up at least at the rate of $1/ \rho$ as $\rho \to 0$.  

\medskip 


We now describe briefly some essential ideas of the proofs of Theorems~\ref{thm1} and \ref{thm1.1}. The starting point is the invariance of the Maxwell equations under a change of variables. More precisely, one has 
\begin{lemma}
   	\label{pre:cha}
   	Let $D, D'$ be two open bounded connected subsets of $\R^3$ and  $F: D \rightarrow D'$ be a bijective, orientation-preserving map such that $F\in C^1(\bar{D}), F^{-1} \in C^1(\bar{D}')$. 
   	Let $\eps, \, \mu \in [L^{\infty}(D)]^{3\times 3}$, and $J \in [L^2(D)]^3$. Assume that $(E, H) \in [H(\curl, D)]^2$ is a solution of  the Maxwell equations
   	\begin{equation}
   	\begin{cases}
   	\nabla \times E = i\omega \mu H & \text{ in } D,\\[6pt]
   	\nabla \times H = -i\omega \eps E +J & \text{ in } D.
   	\end{cases}
   	\end{equation}
   	Set   	$E' = F*E$, $H' = F*H$, $\eps' = F_*\eps$, $\mu' = F_*\mu$, and  $J' = \frac{\nabla F J}{|\det \nabla F|} \circ F^{-1}$  in  $D'$.
   	Then $(E', H') \in [H(\curl, D')]^2$ satisfies 
   	\begin{equation}
   	\begin{cases}
   	\nabla \times E' = i\omega\mu' H' &\text{ in } D',\\[6pt]
   	\nabla \times H' = -i\omega\eps' E' +J' & \text{ in } D'.
   	\end{cases}
   	\end{equation}
	 \end{lemma}

Set 
\begin{equation}\label{def-EHrho}
  (\pbE_{\rho}, \pbH_{\rho}) = (F_{\rho}^{-1}*E_c, \, F_{\rho}^{-1}*  H_c) \quad  \mbox{ in } \R^3. 
\end{equation}
It follows from  Lemma ~\ref{pre:cha} that 
 $(\pbE_\rho, \pbH_\rho) \in [H_{\loc}(\curl, \R^3)]^2$ is the unique (Silver-M\"uller) radiating solution to
   \begin{align}
   \label{equ:aux1}
   \begin{cases}
   \nabla \times \pbE_{\rho} = i\omega \mu_{\rho} \pbH_{\rho}  &\text{ in } \mathbb{R}^3,\\[6pt]
   \nabla \times \pbH_{\rho} = -i\omega \eps_{\rho} \pbE_{\rho} + J_\rho &\text{ in } \mathbb{R}^3,
   \end{cases}
   \end{align}
   where
\begin{equation}\label{epsmurho}
   \big(\eps_{\rho}, \mu_{\rho} \big) = \big( { F_\rho^{-1}}_*\eps_c, {F_\rho^{-1}}_*\mu_c \big) = \left\{\begin{array}{cl} \big(I, I \big) &\text{ in } \mathbb{R}^3\setminus B_{\rho},\\[6pt]
  \big( \rho^{-1}\eps(\cdot / \rho), \rho^{-1} \mu (\cdot/ \rho) \big) &\text{ in } B_{\rho},\end{array}\right.
\end{equation}
   and 
\begin{equation}\label{Jrho}
 J_\rho =  \left\{ \begin{array}{cl}J_{\ext} &\text{ in } \R^3\setminus B_2,\\[6pt] 
\dsp  \rho^{-2}J_{\inte}(\cdot/ \rho) &\text{ in } B_{\rho},   \\[6pt]
   0 & \text{ otherwise}. \end{array}\right.
\end{equation}
We can then derive Theorems~\ref{thm1} and  \ref{thm1.1} by studying the difference between $(\pbE_\rho, \pbH_\rho)$ and $(E, H)$ in $\mR^3 \setminus B_1$ and the behavior of $(\rho \pbE_\rho, \rho \pbH_\rho)(\rho \cdot )$ in $B_1$. As in the acoustic case, this is the study of the effect of a small inclusion. 
It is well-known that when the material parameters inside a small inclusion are bounded from below and above by positive constants, then the effect of a small inclusion is small (see, e.g., \cite{VogeliusVolkov, AVV}). Without this assumption, the effect of the small inclusion might not be  small (see, e.g.,  \cite{Kohn, Ng-1}) unless there is an appropriate lossy layer. In our setting, the boundedness assumption is violated (see \eqref{epsmurho}) and no lossy layer is used. Nevertheless, the effect of the small inclusion is still small  due to the special structure in \eqref{epsmurho}.  It is worth noting that  the definitions of resonance and non-resonance, and the condition of compatibility \eqref{passive}  appear very naturally in this context. Indeed, if $(E_c, H_c)$ is bounded in $[H(\curl, B_1)]^2$ (and $J = 0$ in $B_1$) then one can easily check that, up to a subsequence,  $(\rho \pbE_\rho, \rho  \pbH_\rho)(\rho \cdot) = (E_c, H_c)$ converges weakly in $[H(\curl, B_1)]^2$ to $(E_0, H_0)$ which satisfies  the system \eqref{def:resonance}. 

As in the proof for the scalar case, we  introduce  $(E_{1, \rho}, H_{1, \rho}) \in [H_{\loc}(\curl, \R^3 \setminus B_{\rho})]^2$  as the unique    radiating solution to the system
\begin{equation}
\label{equ:aux2}
\begin{cases}
\nabla \times E_{1, \rho} = i\omega H_{1, \rho}  &\text{ in } \mathbb{R}^3\setminus B_{\rho},\\[6pt]
\nabla \times H_{1, \rho} = -i\omega E_{1, \rho} + J_{\ext} & \text{ in } \mathbb{R}^3\setminus B_{\rho},\\[6pt]
E_{1, \rho} \times \nu= 0 &\text{ on } \partial B_{\rho}. 
\end{cases}\end{equation}
We extend $(E_{1, \rho}, H_{1, \rho})$ by $(0, 0)$ in $B_{\rho}$, and still denote this extension by $(E_{1, \rho}, H_{1, \rho})$. Define
\[
(E_{2, \rho}, H_{2, \rho}) := (E_{1, \rho}, H_{1, \rho}) -(E,H)\quad \mbox{ and } \quad  (E_{3, \rho}, H_{3, \rho}): = (\pbE_{\rho}, \pbH_{\rho}) - (E_{1, \rho}, H_{1, \rho})  \quad \mbox{ in } \mR^3.
\]
Then  $(E_{2, \rho}, H_{2, \rho})\in [H_{\loc}(\curl, \R^3 \setminus B_{\rho})]^2$ is the unique  radiating solution to the system
\begin{equation*}
\begin{cases}
\nabla \times E_{2, \rho} = i\omega H_{2, \rho}  &\text{ in } \mathbb{R}^3\setminus B_{\rho},\\[6pt]
\nabla \times H_{2, \rho} = -i\omega E_{2, \rho} & \text{ in } \mathbb{R}^3\setminus B_{\rho},\\[6pt]
E_{2, \rho} \times \nu= - E\times \nu &\text{ on } \partial B_{\rho},
\end{cases}
\end{equation*}
and $(E_{3, \rho}, H_{3, \rho})\in [\bigcap_{R> 1}H(\curl, B_R \setminus \partial B_{\rho})]^2$ is the unique  radiating solution to the system 
\begin{equation}\label{continue3-???}
\begin{cases}
\nabla \times E_{3, \rho} = i\omega \mu_\rho H_{3, \rho}  &\text{ in } \mathbb{R}^3\setminus \partial B_\rho,\\[6pt]
\nabla \times H_{3, \rho} = -i\omega \eps_\rho E_{3, \rho} + J_{\rho} \mathds{1}_{B_{\rho}} & \text{ in } \mathbb{R}^3\setminus \partial B_\rho,\\[6pt]
[E_{3, \rho} \times \nu] = 0,  \, [H_{3, \rho} \times \nu] = -H_{1, \rho}\times \nu |_{\ext}  &\text{ on } \partial B_\rho. 
\end{cases}
\end{equation}

The problem is then to understand the behavior of 
$(E_{1, \rho}, H_{1, \rho})$, $(E_{2, \rho}, H_{2, \rho})$ and $(E_{3, \rho}, H_{3, \rho})$. After making a change of variables $x \to x/ \rho$, one needs to deal with the Maxwell equations with small frequency. Our analysis is variational and based on various recent results on compactness related to  $H (\curl, \cdot)$ (see \cite{Ng-Maxwell-Superlensing}) and on the Helmholtz decomposition (see \cite[Section 3]{MinhLoc} and also \cite{Girault}).  
We just mention some simple facts related to $(E_{1, \rho}, H_{1, \rho})$ to point out some of the subtleties associated with the electromagnetic nature of the problem. A key ingredient in the study of the behavior of the magnetic fields (after scaling and formally taking the limit $\rho \to 0$) is 
 
\begin{lemma}\label{uniquenessstatic} Assume that $\R^3\setminus D$ is simply connected. Let $u \in H_{\loc}(\curl, \mR^3 \setminus D) \cap H_{\loc}(\dive, \mR^3 \setminus D)$ be such that 
\begin{equation*}
\left\{\begin{array}{cl}
\nabla\times u = 0 &\text{ in } \R^3\setminus D,\\[6pt]
\dive u = 0 &\text{ in } \R^3\setminus D,\\[6pt]
u\cdot \nu = 0 &\text{ on } \partial D,
\end{array} \right. 
\end{equation*}
and 
\begin{equation}\label{US-1}
|u(x)|  = O(|x|^{-2})\mbox{ for large $|x|$}.  
\end{equation}
Then $u = 0$ in $\mR^3 \setminus D$.  
\end{lemma}

A key ingredient in the study of the bahavior of the electric fields (after scaling and formally taking the limit $\rho \to 0$) is

 \begin{lemma}\label{uniquenessstaticE} Assume that $\R^3\setminus D$ is simply connected and $u \in H_{\loc}(\curl, \mR^3 \setminus D) \cap H_{\loc}(\dive, \mR^3 \setminus D)$ is such that 
\begin{equation*}
\left\{\begin{array}{cl}
\nabla\times u = 0 &\text{ in } \R^3\setminus D,\\[6pt]
\dive u = 0 &\text{ in } \R^3\setminus D,\\[6pt]
u\times \nu = 0 &\text{ on } \partial D, 
\end{array} \right. \quad \int_{\Gamma_i} u \cdot \nu = 0 \mbox{  for all connected components $\Gamma_i$ of $\partial D$}, 
\end{equation*}
and 
\begin{equation}\label{US-2}
|u(x)| = O(|x|^{-2}) \mbox{ for large $|x|$.} 
\end{equation}
Then $u = 0$ in $\mR^3 \setminus D$.  
\end{lemma}

It is worth noting that in Lemma~\ref{uniquenessstaticE}, even though one imposes two conditions (on the tangential components) on the boundary, instead of one condition (on the normal component) as in Lemma~\ref{uniquenessstatic}, one also requires the integral conditions on the normal components.  This requirement is indeed necessary, and reveals the intrinsic electromagnetic nature of the problem .  The proof of the first lemma is  based on a representation of $u$ of the form $\nabla \varphi$, and the proof of the second lemma is based on a representation of $u$ of the form $\nabla  \times \varphi$ with appropriate decays of $\varphi$ at infinity. These are invariants of the Helmholtz decomposition in an unbounded domain. Using these representations, the proofs are then just based on  standard integration by parts arguments.

\section{Approximate cloaking for electromagnetic waves in the time regime} \label{sec-5}

This section is devoted to approximate cloaking using transformation optics for Maxwell equations in the {\it time} domain. The results discussed here are from the joint work of  Tran and the first author, \cite{MinhLocT}. To our knowledge, this is the first rigourous work on approximate cloaking for electromagnetic waves in the time domain. 

The cloaking device consists of two layers: a transformation based cloak and a fixed lossy layer as in the spirit of \cite{Ng-Vogelius-3, Ng-Vogelius-4}. 
We assume here that the transformation based cloak in $B_2 \setminus B_1$ is given by \eqref{Frho-M}, the lossy layer in $B_1 \setminus B_{1/2}$ is given by the triple $(I, I, 1)$ modelling the permittivity, permeability, and conductivity,  and the cloaked region $B_{1/2}$ is given by the pair $(\eps_O, \mu_O)$ -- two matrix-valued functions characterizing the permittivity and permeability of the cloaked object. The pair $(\eps_O, \mu_O)$  satisfies the same assumptions as in the previous section.   We assume that the medium is homogeneous outside the cloaking device and the cloaked region. In the presence of the cloaked object and the cloaking device, the medium in the whole space $\R^3$ is described by the triple $(\eps_c, \mu_c, \sigma_c)$ given by 
\begin{equation}\label{medium-cloak}
(\eps_c, \mu_c, \sigma_{c}) =  
\left\{\begin{array}{cl}
(I, I , 0) &\mbox{ in } \R^3\setminus B_2,\\[6pt]
({F_{\rho}}_*I, {F_{\rho}}_*I, 0) & \mbox{ in } B_2\setminus B_{1},\\[6pt]
(I , I, 1) & \mbox{ in } B_1\setminus B_{1/2}, \\[6pt]
(\eps_O, \mu_O, 0) & \mbox{ in } B_{1/2}.
\end{array}\right. 
\end{equation}


Let $\mathcal J$ represent a current. We assume that
\begin{equation}\label{J-1}
\mathcal J\in L^1([0, \infty);[L^2(\R^3)]^3) \mbox{ with } \operatorname{supp} \mathcal J \subset [0, T_0]\times (B_{R_0} \setminus B_2), \mbox{ for some } T_0> 0, \,  R_0>2, 
\end{equation}
and 
\begin{equation}\label{J-2}
\dive \mathcal J = 0 \mbox{ in } \mR_+ \times \mR^3. 
\end{equation}
The electromagnetic field  generated by $\mathcal J$, with zero data at time $0$, and in the presence of the cloaking device,  is the unique weak solution $(\EE_c, \HH_c) \in L^{\infty}_{\loc}([0, \infty); [L^2(\R^3)]^6)$ to the system
\begin{equation}
\label{equ:wcloak}
\begin{cases}
\dsp  \eps_c \frac{\partial \EE_c}{\partial t} =  \nabla \times \HH_c -  \mathcal J -  \sigma_c \EE_c &\text{ in } (0, +\infty)\times\mathbb{R}^3,\\[6pt]
\dsp  \mu_c \frac{\partial \HH_c}{\partial t}  = - \nabla \times \EE_c   &\text{ in } (0, +\infty)\times\mathbb{R}^3,\\[6pt]
\EE_c(0, \cdot) = \HH_c(0,\cdot ) = 0 & \text{ in } \R^3.
\end{cases}\end{equation}
In homogeneous space, the field generated by $\mathcal J$, with zero data at time $0$, is the unique weak solution $(\EE, \HH)\in L^{\infty}_{\loc}([0, \infty); [L^2(\R^3)]^6)$ to the system 
\begin{align}
\label{equ:ncloak}
\begin{cases}
\dsp \frac{\partial \EE}{\partial t}  =  \nabla \times \HH - \mathcal J  &\text{ in } (0, +\infty)\times\mathbb{R}^3,\\[6pt]
\dsp  \frac{ \partial \HH}{\partial t} = - \nabla \times \EE  &\text{ in } (0, +\infty)\times\mathbb{R}^3,\\[6pt]
\EE(0, \cdot) = \HH(0, \cdot ) = 0 & \text{ in } \R^3.
\end{cases}
\end{align}

The definition of weak solution,  is as follows 

\begin{definition}\label{def} Let  $\eps,  \, \mu,  \,  \in [L^{\infty}(\R^3)]^{3\times 3}$~, $\sigma\in L^{\infty}(\R^3)$ be such that $\eps$ and $\mu$ are real,  symmetric, and uniformly elliptic in $\mR^3$,  and $\sigma$ is real and nonnegative  in $\mR^3$, 
and  let $ \mathcal J \in L^1_{\loc}([0, \infty); [L^2(\R^3)]^3)$. A pair  $(\EE, \HH)\in L^{\infty}_{\loc}([0, \infty); [L^2(\R^3)]^6)$ is called a weak solution of 
\begin{equation}\label{equ}
\begin{cases}
\dsp  \eps \frac{\partial \EE}{\partial t} = \nabla\times \HH  - \mathcal J -  \sigma \EE &\text{ in } (0, +\infty)\times \R^3,\\[6pt]
\dsp\mu \frac{\partial \HH}{\partial t} = -  \nabla\times \EE  &\text{ in } (0, +\infty)\times \R^3,\\[6pt]
\dsp \EE(0,\cdot ) =  \HH(0, \cdot ) = 0 & \text{ in } \R^3, 
\end{cases}
\end{equation}
iff  
\begin{equation}\label{def-e1}
\begin{cases}
\dsp \frac{d}{dt}\langle \eps\EE(t, \cdot), E\rangle + \langle \sigma \EE(t, \cdot), E\rangle - \langle \HH(t, \cdot), \nabla\times E\rangle = -\langle \mathcal J(t, \cdot) , E\rangle~, \\[8pt]
\dsp \frac{d}{dt}\langle \mu\HH(t, \cdot ), H\rangle  + \langle \EE(t, \cdot), \nabla\times H\rangle = 0~,
\end{cases} \quad \mbox{ for } t > 0, 
\end{equation}
for all $(E, H)\in [H(\curl, \mR^3)]^2$, and 
\begin{equation}\label{def-e2}
\EE(0, \cdot ) = \HH(0, \cdot ) = 0 \mbox{ in } \R^3.
\end{equation}
\end{definition}

Some comments about Definition~\ref{def} are in order. System \eqref{def-e1} is understood in the distributional sense (in $t$). The initial conditions   \eqref{def-e2}  are understood as 
\begin{equation}\label{trace-sense}
\langle \eps \EE(0, \cdot), E \rangle  = \langle \mu \HH(0, \cdot ), H \rangle = 0  \quad \mbox{ for all } (E, H)\in [H(\curl, \mR^3)]^2. 
\end{equation}
From \eqref{def-e1}, one can check that  
\begin{equation}
\langle \eps \EE(t, \cdot), E \rangle, \langle \mu \HH(t, \cdot), H \rangle \in W^{1, 1}_{\loc}([0, + \infty)). 
\end{equation}
This in turn ensures the existence of the trace  in \eqref{trace-sense}. 

\medskip 
Concerning the well-posedness of \eqref{equ}, we have, see,  e.g., \cite[Theorem 3.1]{Ng-Vinoles}. 
  
\begin{proposition}\label{prop-wellposed} Let $ \mathcal J \in L^1_{\loc}([0, \infty); [L^2(\R^3)]^3)$. There exists a unique weak solution $(\EE, \HH) \in L^{\infty}_{\loc}([0, \infty); [L^2(\R^3)]^6)$ of \eqref{equ}. Moreover, for each $T> 0$, the following estimate holds
	\begin{equation}\label{es}
	\int_{\R^3}|\EE (t, x) |^2 + |\HH (t, x) |^2 dx \leq C\left(\int\limits_{0}^t\Big\| \mathcal J(s,\cdot) \Big\|_{L^2(\R^3)} ds\right)^2 \quad \mbox{ for } t\in [0,T],
	\end{equation}
	for some positive constant $C$ depending only on $T$, and  the ellipticity of  $\eps$ and $\mu$. 
	
\end{proposition}

\begin{remark} \rm  In \cite{Ng-Vinoles}, the authors also considered dispersive materials and hence dealt with Maxwell equations which are non-local in time as well. 
\end{remark}

The precise sense of cloaking is given by 

\begin{theorem}\cite[Theorem 1.1]{MinhLocT}\label{thm-main-1} 
	Let $\rho \in (0, 1)$ and  let $(\EE_c, \HH_c), (\EE, \HH)\in L_{\loc}^{\infty}([0, \infty), [L^2(\R^3)]^6)$ be the unique solutions to systems (\ref{equ:wcloak}) and (\ref{equ:ncloak}) respectively. Suppose \eqref{J-1} and \eqref{J-2} hold, and extend ${\mathcal J}$ by 0 for $t < 0$.  Then, for $K\subset\subset \R^3\backslash \bar B_1$, 
	\begin{equation}\label{main-1}
	\| (\EE_c, \HH_c) - (\EE, \HH)\|_{L^{\infty}((0, T); L^2(K))}\leq C\rho^3\|\mathcal J\|_{H^{11}(\mR; [L^2(\R^3)]^3)},
	\end{equation}	
for some positive constant $C$ depending only on $K$, $R_0$, and $T$.
\end{theorem}

\begin{remark} \rm The order $\rho^3$ in assertion \eqref{main-1} is optimal, and it gives the same degree of visibility as in the  frequency domain. 
\end{remark}

\begin{remark} \rm Estimate \eqref{main-1} requires that $\mathcal J$ be regular. The  degree of regularity of $\mathcal J$ required here might not be optimal. 
\end{remark}

The analysis proceeds in the same spirit as for the acoustic wave equation discussed previously.  We first transform the Maxwell equations in time domain  into a family of Maxwell equations in the time harmonic regime, by taking the Fourier transform of the solution with respect to time. After obtaining appropriate estimates on the  near-invisibility  for the Maxwell equations in the time harmonic regime, we simply invert the Fourier transform. This idea has its roots in the work of Nguyen and Vogelius in \cite{Ng-Vogelius-3} (see also \cite{Ng-Vogelius-4}) in the cloaking context; it was also used to establish the validity of  impedance boundary conditions in the time domain in \cite{MinhLinh}, and recently to study cloaking via a change of variables for the heat equation \cite{Minh-Tu}. To implement this idea, the central issue is to obtain a degree of visibility in which  the dependence on frequency is {\it explicit} and well-controlled.  The analysis involves a variational formulation, a multiplier technique,  and a duality method in some range of the frequency. 
An intriguing fact about the Maxwell equations in the time harmonic regime worthy of mentioning is that the multiplier technique, which plays an important role in the acoustic setting,  does not adapt well in the  very large frequency regime. A duality method  is used instead. 
Another key technical difficulty is to establish the radiation condition for the Fourier transform (in time) of weak solutions to the general Maxwell equations. The resolution of this difficulty  is of interest in itself.

\end{document}